\documentclass[12pt]{article}

\usepackage{hyperref}
\usepackage{amsfonts}
\def\Bbb#1{\mathbb#1}
\newtheorem{theorem}{Theorem}

\newtheorem{lemma}[theorem]{Lemma}
\newtheorem{cor}[theorem]{Corollary}
\newtheorem{rmrk}[theorem]{Remark}

\newcommand{\po}{{\hspace*{-1ex}}{\bf .  }}

\newcommand{\ra}{\rangle}
\newcommand{\la}{\langle}

\newcommand{\R}{\Bbb{R}}

\newcommand{\C}{\Bbb{C\,}}

\newcommand{\proof}{ \noindent{\em Proof:} ~ }
\newcommand{\x}{\noindent}
\newcommand{\vs}{\vspace{.25in}}
\newcommand{\vsss}{\vspace{.1in}}

\newcommand{\re}{\mbox{Re }}

\newcommand{\spa}{\mbox{span}}

\def\im{{\rm Im \,}}
\def\rk{{\rm rank \,}}
\def\ker{{\rm Ker \,}}

\def\be{\begin{equation} }
\def\ee{\end{equation} }
\def\bea{\begin{eqnarray*} }
\def\eea{\end{eqnarray*} }

\def\qed{\ifhmode\unskip\nobreak\fi\ifmmode\ifinner\else\hskip5 pt
\fi\fi\hbox{\hskip5 pt \vrule width4 pt height6 pt depth1.5 pt
\hskip 1pt }}

\title{ Complete real K\"ahler submanifolds 
\\in codimension two}
\author{Luis A. Florit\thanks{IMPA: Estrada Dona Castorina 110,
22460--320, Rio de Janeiro, Brazil; e-mail: luis@impa.br. 
Research partially supported by CNPq.},  \ \ 
 and \,  Fangyang Zheng\thanks{
The Ohio State University: Columbus, OH 43210, USA and IMS, 
Zhejiang University, Hanzhou, China; e-mail:
 zheng@math.ohio-state.edu.
Research partially supported by IHES and a NSF grant.}}
\date{}

\begin{document}
\maketitle

\vspace{-1ex}
\begin{abstract}

Minimal isometric immersions $f: M^{2n}\rightarrow \R^{2n+2}$ 
in codimension two from a complete K\"ahler manifold into Euclidean 
space had been classified in \cite{dg4} for $n\geq 3$. 
In this note we describe the non--minimal situation showing 
that, if $f$ is real analytic but not everywhere minimal, 
then $f$ is a cylinder over a real K\"ahler surface 
$g: N^4\rightarrow \R^6$, that is, $M^{2n}=N^4\times\C^{n-2}$ and
$f=g\times {id}$ split, where ${id} \colon\C^{n-2}\cong\R^{2n-4}$ 
is the identity map. Moreover, $g$ can be further described.

\end{abstract}
\medskip

\vsss
\x {\bf \S 1.  Introduction}
\vs

By a {\em real K\"ahler Euclidean submanifold} we mean a
smooth ($C^\infty$) isometric immersion 
$f: M^{2n} \rightarrow \R^{2n+p}$ 
from a K\"ahler manifold $M^{2n}$ of real dimension $2n$ into 
Euclidean space.
As expected, the K\"ahler structure imposes strong restrictions
on the immersion. In fact, the hypersurface situation ($p=1$) 
is well understood, both locally and globally. Locally, 
by means of an explicit parametrization (\cite{dg1}),
while if $M^{2n}$ is assumed to be complete, we showed in
\cite{fz2} that $f$ must be a {\it cylinder} over a complete 
orientable surface $g\colon N^2\to\R^3$, that is, 
$M^{2n}=N^2\times\C^{n-1}$ and $f=g\times {id}$
split, where \mbox{${id}\colon\C^{n-1}\cong\R^{2n-2}$}
is the identity map.

In codimension $p=2$ the problem becomes far more interesting.
Although few results were known until now in the local case
unless the immersion has rank at most two (\cite{df0}), the complete 
case for dimension $n\geq 3$
is well understood for {\it minimal} immersions. Here, $f$ is 
either a holomorphic complex hypersurface under an identification 
$\R^{2n+2}\cong \C^{n+1}$, 
or a cylinder over a complete minimal real K\"ahler surface 
$g\colon N^4\to\R^6$, or it is essentially completely
holomorphically ruled, i.e., 
$M^{2n}$ is the total space of a holomorphic vector bundle over 
a Riemann surface, and $f$ maps each fiber onto a linear 
subvariety in $\R^{2n+2}$ (\cite{dr}). The paper \cite{dg4} 
was devoted to give a precise description of the latter case  
in terms of a Weierstrass-type representation.
Interesting explicit examples were then given. 

Our main purpose here is to understand the general (analytic) situation, 
that is, to drop the minimality assumption on $f$ and hence to complete 
the global classification by showing that, if non--minimal, $f$ must be 
a cylinder over real K\"ahler surface in $\R^6$:

\begin{theorem}\po\label{cod2} Let $f\colon M^{2n} \rightarrow\R^{2n+2}$, 
$n\geq 2$, be a complete analytic real K\"ahler Euclidean submanifold. 
If $f$ is not everywhere minimal, then 
$M^{2n}=N^4\times \C^{n-2}$ and $f=g\times {id}$ split, for some 
complete real K\"ahler Euclidean surface $g\colon N^4 \rightarrow \R^6$.
\end{theorem}

Although we believe it is superfluous, the analyticity 
assumption appears since we do not know at this point how to 
deal with the possible gluing phenomenon. In general,
we are able to show that there is an open subset 
${\cal U}\subset M^{2n}$
such that ${\cal U}=V^4\times\C^{n-2}$ and the restriction 
$f|_{\cal U}=g\times{id}$ is a cylinder over some 
$g\colon V^4\to\R^6$.

By Hartman's Theorem (\cite{h}), a complete flat Euclidean submanifold 
in codimension two is a cylinder over a flat surface 
$g\colon N^2\to\R^4$, although the decomposition 
\mbox{$f=g\times{id}\colon M^{2n}=N^2\times \R^{2n-2}\to \R^{2n+2}$}
does not need to be a K\"ahler one. The local case has been classified in 
\cite{df1}.

\vsss

The surface case $n=2$ in codimension two is not well 
understood if the immersion is minimal. It was shown in 
\cite{dg3} that there is a family of irreducible 
(i.e., neither a product of 
two surfaces in $\R^3$ nor a cylinder over a surface in $\R^4$) 
complete minimal real K\"ahler surfaces in $\R^6$ which are 
neither holomorphic nor complex ruled. 
However, we can describe the non--minimal situation:

\begin{theorem}\po\label{surface}
Let $f\colon M^4 \to \R^6$ be a complete irreducible analytic real 
K\"ahler Euclidean surface that is not everywhere minimal. Then,  
$M^4=N^2\times \C$, and there is an open dense subset 
$W\subset M^4$ such that, along each connected component $V$ 
of $W$, the restriction $f|_V$ is a composition of analytic 
isometric immersions. That is, $V$ is a real K\"ahler Euclidean 
hypersurface, $g\colon V\to \R^5$, and $f|_V=h\circ g$, where 
$h\colon\, {\cal U}\subseteq \R^5\to\R^6$ is a flat hypersurface,
for some open subset ${\cal U}\subseteq \R^5$ with 
$g(V)\subset {\cal U}$.
\end{theorem}  

In some situations, we can assure that the whole $f$ is
globally a composition; see Remark~\ref{compos}.
Moreover, since the hypersurface situation is parametrically
understood, $f|_V$ above can now also be parametrically described.
More importantly, the proof of Theorem~\ref{surface}, being local
in nature, contains
the ingredients that allow to give a {\it local} parametric
classification of the (not necessarily analytic) everywhere 
non--minimal real K\"ahler Euclidean submanifolds in codimension 
two; cf. Remark \ref{local}. 

\vsss

Theorem \ref{cod2} is based on the main result in this paper
that holds for any (not necessarily analytic)
complete real K\"ahler Euclidean submanifold, in any codimension:

\vsss
\centerline {\it The complex relative nullity foliation is always holomorphic.}
\vsss

\noindent The complex relative nullity distribution $D$ is just 
the maximal complex spaces contained in the relative nullity 
$\Delta$ of the immersion; that is,  $D=\Delta \cap J\Delta$, 
where $J$ stands for the K\"ahler structure of $M^{2n}$ and
$\Delta$ for the kernel of the second fundamental form.
It is easy to check that, on the open subset where $D$
attains its minimal dimension, $D$ is an integrable distribution
with complete totally geodesic leaves in both $M^{2n}$ and the 
Euclidean ambient space. The holomorphicity of $D$ then
imposes strong restrictions on the immersion that, 
in codimension two, allow us to easily conclude Theorem~\ref{cod2}.

\vs
\vsss
\x {\bf \S 2.  The complex relative nullity foliation}
\vs

In this section, we will discuss some general properties about the 
complex relative nullity foliation of any complete real K\"ahler
Euclidean submanifold. The main result is that the foliation is always 
a holomorphic one. We then apply this property to the codimension two
case to easily get Theorem \ref{cod2}. Hopefully the main result 
might be useful in other contexts as well.

\vsss
Let $f: M^{2n} \rightarrow \R^{2n+p}$ be a complete real K\"ahler
Euclidean submanifold. We will denote by $\Delta_x$ the {\it relative
nullity} of $f$ at $x\in M$, that is, the nullity space of the second
fundamental form $\alpha(x)$ of $f$ at $x$,
$\Delta_x=\{Z\in T_xM:\alpha(Z,T_xM)=0\}$, 
and by $\nu (x)$ the {\it index of relative nullity} of $f$ at $x$,
i.e., $\nu (x) = \dim_\R\Delta_x$.
Let $U\subset M$ be the open set where $\nu$ attains its minimum value
$\nu_0$,
$$U=\{ x\in M \mid \nu (x) = \nu_0 \}.$$
It is a well-known fact (see e.g. \cite{fe}) that 
$\Delta$ is smooth and integrable in $U$
with totally geodesic leaves in both $M^{2n}$ and
$\R^{2n+p}$. If, in addition, $M$ is complete,
then the leaves are also complete and thus each one is mapped by $f$ 
onto a linear subvariety (that is, a translation of linear subspace)
of $\R^{2n+p}$. 

Let us define the complex subspaces $D_x\subset T_xM$ by
$$D_x = \Delta_x \cap J \Delta_x,$$
where $J$ is the almost complex structure of $M$, and by 
$\nu'(x)=\dim_{\C}D_x$ its complex dimension. 
Observe that $D=\Delta\cap\Delta_J$ as well, where $\Delta_J$ is the 
{\it pluriharmonic nullity} of~$f$ defined in \cite{fz1} by
$\Delta_J:=\{Z\in TM:\alpha(JZ,Y)=\alpha(Z,JY), \forall\, Y\in TM\}$.

Let $\nu_0'$ be the minimum value of $\nu'(x)$ for all $x\in U$, 
and $U_0\subseteq U$ be the open subset where $\nu'=\nu_0'$.
Clearly, since $J$ is parallel, $D$ is also smooth, integrable and
totally geodesic in $U_0$ and its leaves are again complete.
Therefore, each leaf is isometric
to~$\C^{\nu_0'}$ and is mapped by $f$ onto a linear subvariety in
$\R^{2n+p}$. We will call the leaves of~$D$ in~$U_0$ the
{\em complex relative nullity} foliation of $f$ from now on. 

Let us write $r=n-\nu_0'$ and fix any $x\in U_0$. Let $V$ be the 
space of type $(1,0)$ tangent vectors at $x$, that is, $V$ is the
complex subspace of $(T_xM) \! \otimes \! \C$ defined as 
\mbox{$V=\{v-iJv: v\in T_xM\}$}.
Denote by $W=W_x\cong \C^r$ the complex linear subspace of $V$ 
perpendicular to $D_x$, that is, 
$$W\oplus \overline{W} = D_x^{\perp } \otimes \C.$$
Let $C\colon D\times D^\perp \to D^\perp$ be the {\it twisting tensor}
(also called the {\it splitting tensor}) of 
the totally geodesic foliation $D$ defined by
$C_TX=C(T,X)=-(\nabla_XT)_{D^\perp}$, where $\nabla$ stands for the
Levi-Civita connection of $M$ and $(\ )_{D^\perp}$ for the orthogonal
projection onto $D^\perp$. Fix $T\in D$ and write for the 
complexified operator $C_T$,
$$C_T(e_i) = \sum_{j=1}^r ( A_{ij}e_j 
+ \overline{B_{\overline{i}j}} \overline{e_j} ),$$
for a basis ${\cal B}=\{ e_1, \ldots , e_r\}$ of $W$. 
We need a basic property of the twisting tensor:

\begin{lemma}\po\label{1} 
Let $f: M^{2n} \rightarrow \R^{2n+p}$ be 
a complete real K\"ahler submanifold and $x\in U_0$. Then, for any
$T\in D_x$  and any $\lambda \in \C$, the matrix
$$\tilde{C} = \left[ \begin{array}{ll} 
\lambda A & \overline{\lambda B} \\ 
\lambda B & \overline{\lambda A} 
\end{array} \right]$$
has no non-zero real eigenvalues.
\end{lemma}

\proof Write $\lambda = a+ib$ and take a geodesic $\gamma$ with
$\gamma'(0)=aT+bJT$. From $JD=D$ we have $C_{JT}=JC_T$. So, under the
frame $\{ e_i, \overline{e_i}\}$ for $D^{\perp } \otimes \C$, 
the twisting tensor $C_{\gamma'(0)}$ is represented by the above 
$2r\times 2r$ matrix. Since $D$ is totally
geodesic and is contained in the nullity space of the curvature tensor, 
it is easy to check that $C=C_{\gamma'(t)}$ satisfies the Riccati 
equation $C'=C^2$. Hence, since the leaves of $D$ are complete, 
$C_{\gamma'(0)}$ cannot have any non-zero real eigenvalue, just like 
the case for conullity operators as observed by Abe (cf. \cite{a}):
the solution of the above Riccati equation is 
$C(t)=C(0)(I-tC(0))^{-1}$.
\qed

\vs
Let us now recall the following decomposition of the second 
fundamental form $\alpha$ of~$f$ at $x\in M$ (see \cite{fhz}). 
Extend $\alpha$ bilinearly over $\C $, and still denote it
by $\alpha$,
$$\alpha : (T_xM)\! \otimes \! \C  \times (T_xM)\! \otimes \! 
\C  \rightarrow T^\perp_xM\! \otimes \! \C.$$
 Using that
$(T_xM)\! \otimes \! \C = V \oplus \overline {V}$, we can write
$$H=\alpha|_{V\times\overline{V}} \ \ \ \ \ \mbox{and} \ \ \ \ \ \
S=\alpha|_{V\times V}$$
for the $(1,1)$ and $(2,0)$ parts of $\alpha$, respectively. 
Let $W'\subset V$ be the complex linear subspace given by
$$W' = \mbox{ker}\,H \ \cap \ \mbox{ker}\,S.$$
Hence, $D\otimes \C = W' \oplus \overline{W'}$ and $V=W\oplus W'$.  
With respect to the basis ${\cal B}$, we have that $A$ and $B$ in
Lemma~\ref{1} are $r\times r$ complex
matrices, while $H$ and $S$ are Hermitian and complex symmetric
matrices, respectively, with values in the (complexification of the)
normal space of $f$ at $x$. These operators satisfy the following
compatibility conditions:

\begin{lemma}\po\label{2} 
At any $x\in U_0$ and under any basis ${\cal B}$ of $W$, the matrices
$AS$ and $BH$ are always symmetric. Moreover, it holds that $AH=SB^t$.
\end{lemma}

\proof Note that under a base change, say 
${e_i'} = \sum_{j=1}^r P_{ij}e_j$, 
for $P\in GL(r,\C)$, the matrices $A$, $B$, $H$, $S$ change to 
\be\label{eq1}
PAP^{-1},\ \ \ \overline{P}BP^{-1},\ \ \ PHP^{\ast },\ \ \ PSP^t,
\ee
respectively. So the symmetry of $AS$ and $BH$, as well as 
the identity $AH=SB^t$, are independent of the choice of the frame. 

Let $\{ e_i, e_{\alpha }\}$ be a local unitary frame of type $(1,0)$
tangent vector fields near $x$, such that $\{ e_i\}$ gives a basis of
$W$ in a neighborhood of $x$. Here and below we will use
the convention that Latin indices $i$, $j$,... will run from 
$1$ to $r$, while Greek indices $\alpha, \beta$,... will 
run between $r+1$ and $n$. Also, let $a$, $b$,... 
run through the full range, from $1$~to~$n$. 

Let $g$ and $\theta$ be the matrices of metric and connection 
of $M$ under the frame $\cal B$. Also, write
$$\xi_{a\mu } = \langle\tilde \nabla e_a , w_{\mu } \rangle,$$
where $1\leq \mu \leq p$, $\{ w_1, \ldots , w_p\}$ is an
orthonormal frame of the normal bundle of $f$ near~$x$, and 
$\tilde \nabla$ is the covariant differentiation in $\R^{2n+p}$. 
we get
$$\xi_{\alpha \mu }=0, \ \ \ 
\xi_{i\mu } = \sum_{k=1}^r (S_{ik}^{\mu } \varphi_k +
H_{i\overline{k}}^{\mu } \overline{\varphi_k}),$$
where $\{ \varphi_a\} $ is the dual coframe of $\{ e_a\}$. 
Fix any $\alpha$, write 
$e_{\alpha } = \frac{1}{2}(T_\alpha-\sqrt{-1}JT_\alpha)$ and
$C_{T_\alpha}(e_i)=\sum_{j=1}^r (A^\alpha_{ij}e_j+ \overline{B^\alpha_{\overline{i}j}}\,\overline{e_j})$. 
 From now on, for the sake of simplicity, we omit  
the superscript $\mu$ for $S$ and $H$. Then we obtain
$$\theta_{\alpha i} = \la\nabla e_\alpha,e_i\ra=- 
\sum_{j=1}^r (A^\alpha_{ji} \varphi_j -B^\alpha_{\overline{j}i} \overline{\varphi_j}).$$
By the Codazzi equation, we get
\begin{eqnarray*}
 0 & = & - d\xi_{\alpha\mu} \ = \ 
- \sum_{i=1}^r (\theta_{\alpha i}\, \xi_{i\mu}) \\
 & = &  \sum_{i,j,k=1}^r (A^\alpha_{ji}\,\varphi_j + B^\alpha_{\overline{j}i}\,\overline{\varphi_j})\wedge 
(S_{ik}\,\varphi_k + H_{i\overline{k}}\,\overline{\varphi_k} ) \\
& = &\sum_{j,k=1}^r\left((A^\alpha S)_{jk}\, \varphi_j\wedge \varphi_k
+ (B^\alpha H)_{\overline{j}\overline{k}}\,
\overline{\varphi_j} \wedge \overline{\varphi_k} +
(A^\alpha H-S(B^{\alpha})^t)_{j\overline{k}}\, \varphi_j \wedge \overline{\varphi_k}\right).  
\end{eqnarray*}
We conclude that $A^\alpha S$ and $B^\alpha H$ are both symmetric, 
and $A^\alpha H= S(B^\alpha)^t$ for all $\alpha$. 
\qed

\vs

Next, we observe that the operators $A$, $B$ also satisfy a compatibility 
condition with the curvature tensor $R$ of $M$. 

\begin{lemma}\po\label{3}For any 
$1\leq i,j,k,l\leq r$, the components $A$ and $B$ of $C_T$ satisfy
that
$$\sum_{p=1}^r A_{ip} R_{p\overline{j}k\overline{l}}=
\sum_{p=1}^r A_{kp} R_{p\overline{j}i\overline{l}},
\ \ \ \ \mbox{and} \ \ \ \ \ \ 
\sum_{p=1}^r B_{\overline{i}p} R_{p\overline{j}k\overline{l}} =
\sum_{p=1}^rB_{\overline{j}p} R_{p\overline{i}k\overline{l}}.$$
\end{lemma}

\proof 
Recall that, since $M$ is K\"ahler, 
$R(V,V)=R(\overline V,\overline V)=0$. Moreover, the relative
nullity, and hence $D$, is always contained in the nullity of $R$ 
by the Gauss equation. So, the second Bianchi identity gives us that
$0=(\nabla_{e_i}R)(T,\overline e_j, e_k, \overline e_l) -
(\nabla_{e_k}R)(T,\overline e_j, e_i, \overline e_l) =
R(C_Te_i,\overline e_j, e_k, \overline e_l)-
R(C_Te_k,\overline e_j, e_i, \overline e_l),$
which is the first relation we wanted to prove. 
The proof of the second one is similar.
\qed

\vs

Note that Lemma \ref{3} holds true under any basis ${\cal B}$ of $W$,
not necessarily an orthogonal one. The symmetry in Lemma \ref{3} was
observed in \cite{wz}, where the situation is intrinsic. 
In our case here, 
the manifold $M$ may not have nonpositive or nonnegative bisectional
curvature, and the complex relative nullity $D$, although is always
contained in the nullity of $M$, may not coincide with the
nullity. So, in order for us to exploit techniques of the proof of 
Theorem A in \cite{wz}, we need more symmetry conditions on the
components $A$ and $B$ of the twisting tensor.

By the Gauss equation, the curvature tensor $R$ of $M$ is given by
$$\la R(X,\overline{Y})Z,\overline{U}\ra =
\langle H{(X,\overline{U})}, H({Z,\overline{Y}}) \rangle 
- \langle S({X,Z}) , \overline{S({Y,U})} \rangle,$$
for all $X,Y,Z,U\in V$. Let us introduce the tensor $\hat R$ by
$$\la \hat R(X,\overline{Y})Z,\overline{U}\ra
=\langle H{(X,\overline{U})}, H({Z,\overline{Y}}) \rangle 
+ \langle S({X,Z}) , \overline{S({Y,U})} \rangle,$$
also with $\hat R(V,V)=\hat R(\overline V,\overline V)=0$.
It has all the symmetries of $R$, i.e., it is~a~curvature-like
tensor. Taking a unitary basis $\{w_1,,\dots,w_r\}$ of $W$,
the Ricci tensor of $M$ is given~by
$$Q(X,Y) = \mbox{Ric}(X,\overline Y) = 
\sum_{j=1}^r \la R(X,\overline w_i)w_i,\overline Y\ra=
\sum_{j=1}^r \la R(X,\overline Y)w_i,\overline w_i\ra= (Q^H-Q^S)(X,Y),
$$
with
$Q^H(X,Y)=\sum_{i=1}^r\la H(X,\overline w_i),\overline{H(Y,\overline w_i)}\ra$
and
$Q^S(X,Y)=\sum_{i=1}^r\la S(X, w_i), \overline{S(Y,w_i)}\ra$.
Notice that the corresponding $\hat{Q} = Q^H+ Q^S $ for $\hat R$ 
is positive definite on $W$, 
since both $Q^H$ and $Q^S$ are positive semidefinite, and $W$ is 
the orthogonal complement of the common nullity of $H$ and $S$. 

Since $AS$ and $BH$ are symmetric by Lemma \ref{2}, Lemma \ref{3}
still holds true if we replace $R$ by $\hat{R}$. That is, we have:
\begin{lemma}\po\label{4}  
For any $1\leq i,j,k,l\leq r$, the components $A$, $B$ 
of the twisting tensor satisfy
$$\sum_{p=1}^r A_{ip} \hat{R}_{p\overline{j}k\overline{l}} = 
\sum_{p=1}^r A_{kp} \hat{R}_{p\overline{j}i\overline{l}}\,, 
\ \ \ \ \mbox{and} \ \ \ \ \
\sum_{p=1}^r B_{\overline{i}p} \hat{R}_{p\overline{j}k\overline{l}} = 
\sum_{p=1}^r B_{\overline{j}p} \hat{R}_{p\overline{i}k\overline{l}}.$$ 
In particular, the matrix $B\hat{Q}$ is always symmetric.
\end{lemma}

The last sentence is the result of contracting the second identity 
by the inverse of the metric, $g^{\overline{l}k}$. We have all the
ingredients conclude of our main result:

\begin{theorem}\po\label{10} For any complete real K\"ahler
Euclidean submanifold $f: M^{2n}\rightarrow \R^{2n+p}$, 
the complex relative nullity $D=\Delta \cap J\Delta$ in $U_0$ 
is a holomorphic foliation, that is, $B=0$.
\end{theorem}

\proof Consider a basis ${\cal B}$ such that 
$\hat Q(e_i,e_j)=\delta_{ij}$. Under ${\cal B}$,
by the last assertion in Lemma \ref{4}, $B$ becomes
a complex symmetric matrix. By (\ref{eq1}), it is well known 
that there is a basis ${\cal B}'$ under which $B$ is diagonal 
with nonnegative diagonal entries.
The proof now ends just as the one for Theorem A in \cite{wz},
in view of Lemma \ref{1} and Lemma \ref{2}.
\qed
\vs

Let us put together what we know in general about the
complex relative nullity distribution. For this, denote by
${\rm Im}\,H = {\rm span}\,H = 
{\rm span} \{H(X,\overline Y): X,Y\in V\}$.

\begin{cor}\po\label{fund} In the situation of Theorem \ref{10}, 
for any $T\in D$, $A$ is nilpotent, $AS$ is symmetric, and $AH=0$;
that is, the complexified twisting tensor $C_T$ satisfies that 
\be\label{eqf1}
\alpha(C_TX,Y)=\alpha(X,C_TY) \in ({\rm Im}\,H)^\perp, 
\ \ \ \forall\, X,Y\in W,
\ee
\be\label{eqf2}
C_T\colon W\to W, \ \ \ \ C_T^r=0, \ \ \ \ 
{\rm Im}\, C_T \subset {\rm ker}\,H.
\ee
In other words, for the {\rm real} twisting tensor 
$C_T\colon D^\perp \to D^\perp$ we~have that
\be\label{eqf3}
\alpha(C_TX,Y)=\alpha(X,C_TY), \ \ \ \forall\, X,Y\in TM,
\ee
\be\label{eqf4}
C_T\circ J = J \circ C_T,\ \ \ \ C_T^r=0, \ \ \ \ 
{\rm Im}\, C_T \subset \Delta_J.
\ee
\end{cor}
\proof It follows from Theorem \ref{10}, Lemma \ref{1},
Lemma \ref{2}, and the relation 
\be\label{2fhz}
\la S_{ij},H_{k\overline s}\ra=\la S_{kj},H_{i\overline s}\ra
\ \ \ \ \forall\ i,j,k,s, 
\ee
that is an easy consequence of the curvature 
symmetries; see (2) in \cite{fhz}.
\qed

\begin{cor}\po\label{maincor} 
With the hypothesis of Theorem \ref{10}, assume 
that one of the following holds:
$$(a)\  {\rm ker}\,H\subseteq{\rm ker}\,S, 
\ \ \ \ \ \mbox{or} \ \ \ \ \
(b)\ {\rm dim\ Im}\, H \geq p-1.$$
Then, $C=0$. That is, each connected component $U_i$ 
of $U_0$ is isometric to a product 
$$U_i=N^{2r}\times \C^{n-r}, \ \ \ \ r\leq p,$$
and $f|_{U_i} = f'\times id$
split, where $f'\colon N^{2r}\to \R^{2r+p}$ is a real 
K\"ahler Euclidean submanifold, and 
$id\colon \C^{n-r}\cong \R^{2(n-r)}$ is the identity map.
\end{cor}
\proof 
If $(a)$ holds, then ker $H=W'$ and by (\ref{eqf2}) $C$ must vanish.
Assume that $(b)$ holds but $(a)$ does not. By (\ref{2fhz}) we have
that 
\be\label{kernel}
S({\rm ker}\,H,V)\subset ({\rm Im}\, H)^\perp.
\ee
Hence, dim $({\rm Im}\, H)^\perp = 1$, and there is $Z\in {\rm ker}\,H$
such that $S(Z,V)= ({\rm Im}\,H)^\perp$. Consider the hyperplane
$L={\rm ker}\,S(Z,\,\cdot\,)\subset V$. Since from the curvature 
symmetries the relation 
$\langle S_{ij}, S_{ks} \rangle = \langle S_{kj}, S_{is} \rangle$
always holds (see (3) in \cite{fhz}),
we get $S(L,V)\subset {\rm Im}\,H$ and 
\be\label{l}
W' = L \cap {\rm ker\,} H.
\ee
Hence, by (\ref{eqf1}) we obtain that
$\alpha(C_TL,V)\in {\rm Im}\,H \cap ({\rm Im}\,H)^\perp = 0$
since $H$ is Hermitian. This yields $L\subset {\rm ker}\, C$ 
by the last relation in (\ref{eqf2}). Since $C_T$ is nilpotent, 
we conclude again from the last relation in (\ref{eqf2}) that 
${\rm Im}\, C_T \subset L\cap {\rm ker\,} H$. Therefore,
by (\ref{l}), $C=0$.

The estimate on $r$ follows from (\ref{l}) and Lemma 7 of 
\cite{fhz}, where it was proved that 
dim~ker~$H\geq n\, -$~dim~Im~$H$.
\qed

\begin{rmrk}\po {\rm
We point out that, where ker\,$H$ attains its minimal 
possible complex dimension $n-p$,  it holds that 
${\rm ker}\,H\subseteq{\rm ker}\,S$ (see Lemma 7 in \cite{fhz}).}
\end{rmrk}

{\it Proof of Theorem \ref{cod2}:} It follows from Corollary 
\ref{maincor} $(b)$ and the fact that $f$ is minimal if and 
only if $H=0$ (cf. Remark 8 in \cite{fhz}).
\qed

\vs
\x {\bf \S 3. Real K\"ahler surfaces and the local 
case in codimension two}
\vs

Here we basically argue locally to understand real K\"ahler
Euclidean surfaces in $\R^6$, giving the proof of Theorem
\ref{surface}. As a side effect, we 
obtain a local classification of all real K\"ahler 
Euclidean submanifolds in codimension two.
\vsss

{\it Proof of Theorem \ref{surface}:} 
By Corollary \ref{maincor}, $r=2-\nu'_0=2$ and then $D=0$. 
Since $f$ is not minimal, 
we get that $\nu_J=\dim_{\C}\Delta_J=0$ or $1$ almost 
everywhere. In the first case, the composition structure 
of $f$ follows from Theorem 1 in \cite{fz1}, with 
nullity index $\mu=2$ of the curvature tensor of $M^4$ and 
relative nullity $\nu=1$ almost everywhere.
Hence, assume also that $\nu_J=1$ almost everywhere.
The next arguments hold along (connected components of) an 
open dense subset by the analyticity~of~$f$.

Consider $\xi$ a unit (analytic) vector field spanning the line 
bundle $\im H$. Observe that we can choose $\xi$ to be real since 
$H(X,\overline X)$ is real, $X\in V$. 
Take $\{\xi,\eta\}$ an orthonormal basis of the normal space of 
$f$, with corresponding shape operators $A_\xi$ and $A_\eta$. 
By (\ref{2fhz}) we have that 
\be\label{11}
1\leq \rk A_\xi \leq 2,
\ee
and that (\ref{kernel}) holds.
Hence, trace $A_\eta=0$ and $A_\eta\neq 0$ since $D=0$.

The curvature tensor symmetries also imply that
$\la S_{ij},S_{ks}\ra=\la S_{kj},S_{is}\ra$, for all $i,j,k,s$
(cf. (3) in \cite{fhz}). In terms of the complex bilinear 
form $s(X,Y)=\la S(X,Y),\eta\ra$, this is equivalent to
$s(X,Y)s(Z,W) = s(X,W)s(Z,Y)$, for all $X,Y,Z,W\in~V$. 
Taking $Y=X$ with $s(X,X)\neq 0$ and $s(Z,X)=0$, we 
conclude from $\im H \perp \eta$ that
\be\label{12}
\rk A_\eta = 2.
\ee
Since $D=0$, we also get from (\ref{11}) and (\ref{12}) 
that $\im A_\xi\not\subset \im A_\eta$, and, since $\mu$
is even, by the Gauss equation we have that
\be\label{disj}
\im A_\eta \cap \im A_\xi=0.
\ee
The fact that $\mu=0$ is then equivalent, by Gauss equation, 
to $\rk A_\xi = 2$. In this case, we claim that $f$ would split 
as a product of two surfaces in $\R^3$, which is a contradiction.

To prove the claim, consider $\xi_1=\xi, \xi_2=\eta$, and 
take $X,Y\in \ker A_{\xi_i}$. The Codazzi equation for
$A_{\xi_i}$ says that 
\be\label{ker}
A_{\xi_i}[X,Y]=(-1)^i A_{\xi_j}(\psi(X)Y-\psi(Y)X),
\ \ \ \ \ 1\leq i\neq j\leq 2,
\ee
where $\psi(X) = \la \nabla^\perp_X{\xi},\eta\ra$.
 From this, (\ref{12}), (\ref{disj}) and $\rk A_\xi = 2$
we easily get that $\ker A_{\xi_i}$ is integrable and that 
$\xi_i$ is parallel, $i=1,2$. Then, by Ricci equation, we 
obtain that $\ker A_{\xi_i} = \im A_{\xi_j}$. 
On the other hand, Codazzi 
equation for $X\in \ker A_{\xi_i}$ and $Y\in \im A_{\xi_i}$ 
gives $\nabla_X A_{\xi_i}Y = A_{\xi_i}[X,Y]\in \im A_{\xi_i}$.
Therefore, we have the decomposition $TM=\im A_\xi \oplus \im A_\eta$
into orthogonal parallel distributions. The claim now follows from 
the local de Rham decomposition Theorem and the Main Lemma in \cite{m}.

 So, we must have, again, $\mu=2$ and $\nu=1$,
with $\rk A_{\xi}=~1$. In this situation, we obtain the (local)
composition structure of $f$ as follows. 
By (\ref{disj}) and (\ref{ker}) for $\xi_i=\xi$, $\xi_j=\eta$, 
we obtain that
\be\label{ppp}
\ker A_\xi\subset \ker \psi.
\ee
Thus, it is easy to check that $A_\eta$ is a Codazzi tensor, 
and then, since $\rk A_\xi=1$, there is an Euclidean hypersurface 
$g$ whose second fundamental form is $A_\eta$. Now, by (\ref{ppp}), 
we have $(\tilde\nabla_{TM}\xi) \cap \spa\{\eta\}=0$. Therefore, we 
conclude that $f$ is a composition from Proposition 8 in \cite{df2}. 
Observe that the fact that the line bundles spanned by $\xi$ and 
$\eta$ are analytic implies that the immersions $h$ and $g$ also 
are.

It remains only to argue that $M^4 = N^2\times \C$ if $\mu=2$.
It is well known that the nullity $\Gamma$ of the curvature
tensor of any Riemannian manifold is an integrable totally
geodesic distribution, in any open subset where $\mu=\dim \Gamma$
is constant. Moreover, along the (open) set where $\mu$ is 
minimal, in our case $\mu^{-1}(2)$, the leaves are complete.
The twisting tensor $\hat C$ of $\Gamma$
also satisfies the Riccati equation $\hat C_T'=\hat C_T^2$
for any $T\in \Gamma$, since $R(X,T)T=0$ for all $X\in TM$.
Again by the completeness of the leaves of $\Gamma$, this equation 
has no real eigenvalues. But since $\Gamma$ is $J$-invariant,
the same argument that of Lemma \ref{1} gives $\hat C=0$.
The global splitting follows from the analyticity of $f$.
\qed

\begin{rmrk}\po\label{compos}
{\rm We conclude that $f|_W$ in 
Theorem~\ref{surface} is itself a composition if $f$ is an 
embedding. But we do not know if such an $f$ itself is always 
a composition on the whole~$M^4$, $f=h\circ g$. However, 
if $f$ is nowhere flat, then $\rk A_\eta =2$ everywhere. 
So, the normal subbundles spanned by $\xi$ and $\eta$ are 
globally well defined. Hence, if $M^4$ is simply connected,
$g$ is also globally well defined and then it is a cylinder 
over a surface in $\R^3$ by \cite{fz2}. The same holds if $f$ 
is nowhere minimal because $\xi$ would never vanish, but in 
this case, since $\alpha$ decomposes regularly, $h$ is also 
globally well defined, and thus $f=h\circ g$ is globally a 
composition; see Proposition 8 in \cite{df2}.}
\end{rmrk}

\begin{rmrk}\po\label{local}{\rm 
The proof of Theorem \ref{surface} contains all the 
ingredients to obtain the {\it local} classification of nowhere 
minimal real K\"ahler Euclidean submanifolds in codimension 
two $f\colon M^{2n}\to \R^{2n+2}$, 
not necessarily analytic. In terms of the index of 
nullity $\mu$ of $M^{2n}$ and index of relative nullity 
$\nu\leq \mu$ of $f$, the restriction of $f$ to each connected 
component $V$ of an open dense subset $W\subset M^{2n}$
must then be:
\begin{itemize}
\item[$\bullet$] $\mu=2n$: this is the flat case, classified 
parametrically in Theorem 13 in \cite{df1};
\item[$\bullet$] $\mu=\nu=2n-2$: either $f$ reduces
codimension and is then a hypersurface (classified 
parametrically in \cite{dg1}), or is a cylinder over a 
surface in $\R^4$. Otherwise, by Theorems 25 and 27 in 
\cite{df0}, it would be minimal, admitting 
a Weierstrass-type representation; 
\item[$\bullet$] $\mu=2n-2$, $\nu=2n-3$: similarly as in
Theorem \ref{surface}, $f$ is a composition, 
so it reduces to the hypersurface situation;
\item[$\bullet$] $\mu=2n-4$: here $f$ is a product of two
Euclidean hypersurfaces (see Theorem 1 in \cite{fz1}).
\end{itemize}
Observe that $\mu<2n-4$ cannot occur, since the submanifold 
would then be a holomorphic hypersurface 
in $\C^{n+1}\cong \R^{2n+2}$ by \cite{d}, and hence minimal. 
Therefore, the local classification problem of real K\"ahler 
Euclidean submanifolds $f$ in codimension two reduces, aside 
from the gluing phenomena, to the minimal case with $\nu=\mu=2n-4$.
However, we recall that any minimal real K\"ahler submanifold
$f\colon M^{2n}\to \R^{N}$ is the real part of a holomorphic 
complex submanifold $g\colon M^{2n}\to \C^{N}=\R^N\oplus\R^N$, 
$f=\re g$; cf. Theorem~1.11 in \cite{dg1} and Remark~8 in 
\cite{fhz}.}
\end{rmrk}

\end{document}